\documentclass[11pt]{amsart}

\usepackage{amsmath,amssymb}
\usepackage{indentfirst}
\usepackage[all]{xy}
\usepackage{color}
\newtheorem{tm}{Theorem}[section]
\newtheorem{lemma}[tm]{Lemma}

\begin{document}

\renewcommand{\P}{\mathbb{P}}
\newcommand{\R}{\mathbb{R}}
\newcommand{\N}{\mathbb{N}}
\newcommand{\Z}{\mathbb{Z}}
\newcommand{\Q}{\mathbb{Q}}
\newcommand{\scf}{\mathcal{F}}
\newcommand{\scc}{\mathcal{C}}
\newcommand{\scs}{\mathcal{S}}
\newcommand{\scu}{\mathcal{U}}
\newcommand{\scv}{\mathcal{V}}
\newcommand{\ra}{\rightarrow}
\newcommand{\ov}{\overline}
\newcommand{\bs}{\backslash}
\renewcommand{\mid}{\big\vert}
\newcommand{\al}{\alpha}
\newcommand{\da}{\delta}
\newcommand{\e}{\varepsilon}
\newcommand{\et}{\emptyset}
\newcommand{\ga}{\gamma}
\newcommand{\vp}{\varphi}
\newcommand{\Sa}{\Sigma}
\newcommand{\sa}{\sigma}
\newcommand{\la}{\lambda}
\renewcommand{\mod}{\mathrm{mod}}

\newcommand{\mesh}{\mathop{\mathrm{mesh}}}
\newcommand{\ord}{\mathop{\mathrm{ord}}}
\newcommand{\card}{\mathop{\mathrm{card}}}
\newcommand{\diam}{\mathop{\mathrm{diam}}}
\newcommand{\im}{\mathop{\mathrm{im}}}
\newcommand{\st}{\mathop{\mathrm{st}}}
\newcommand{\ANR}{\mathop{\mathrm{ANR}}}
\renewcommand{\int}{\mathop{\mathrm{int}}}
\newcommand{\inv}{^{-1}}
\newcommand{\Sh}{\mathop{\mathrm{Sh}}}
\newcommand{\UV}{\mathop{\mathrm{UV}}}
\newcommand {\Tor}{\mathop{\mathrm{Tor}}}

\title{A proof of the Edwards-Walsh Resolution Theorem without Edwards-Walsh CW-complexes}
\author{Vera Toni\'{c}}
\address{Department of Computer Science and Mathematics\\
Nipissing University\\
100 College Drive, Box 5002\\
North Bay, Ontario P1B 8L7\\
Canada}
\email{vera.tonic@gmail.com}

\date{16 July 2012}

\keywords{Bockstein basis,  cell-like map, cohomological dimension, CW-complex, dimension, Edwards-Walsh resolution, inverse sequence, simplicial complex}

\subjclass[2010]{Primary \textbf{54F45, 55M10}, 55P20, 54C20}

\maketitle \markboth{V.\ Toni\'{c}}{A proof of the Edwards-Walsh Resolution Theorem without Edwards-Walsh CW-complexes}

\begin{abstract}
In the paper  titled ``Bockstein basis and resolution theorems in
extension theory'' (\cite{To}),
we stated a theorem that we claimed to be a generalization of the Edwards-Walsh resolution theorem.
The goal of this note is to show that the main theorem from \cite{To} is in fact equivalent to the Edwards-Walsh resolution theorem, and also that it can be proven without using Edwards-Walsh complexes. We conclude that the Edwards-Walsh resolution theorem can be proven without using Edwards-Walsh complexes.
\end{abstract}

\section{Introduction}

In the paper  titled ``Bockstein basis and resolution theorems in
extension theory'' (\cite{To}),
the following theorem is proven.

\begin{tm} \label{T} Let $G$ be an abelian group with
$P_G=\mathbb{P}$, where $P_G=\{ p \in \P: \Z_{(p)}\in$ Bockstein Basis $ \sigma(G)\}$. Let $n\in \N$ and let $K$ be a connected \emph{CW}-complex with $\pi_n(K)\cong G$, $\pi_k(K)\cong 0$ for $0\leq k< n$. Then for every compact metrizable
space $X$ with $X\tau K$ (i.e., with $K$ an absolute extensor for
$X$), there exists a compact metrizable space $Z$ and a surjective
map $\pi: Z \rightarrow X$ such that
\begin{enumerate}
\item[(a)] $\pi$ is cell-like,
\item[(b)] $\dim Z\leq n$, and
\item[(c)] $Z\tau K$.
\end{enumerate}
\end{tm}


This theorem turns out to be equivalent to the Edwards-Walsh resolution theorem, first stated by R.\ Edwards in \cite{Ed}, with proof published by J.\ Walsh in \cite{Wa}:

\begin{tm}\emph{(R.~Edwards - J.~Walsh, 1981)} \label{EdWa}
For every compact metrizable space $X$ with $\dim_{\Z} X
\leq n$, there exists a compact metrizable space $Z$ and a
surjective map $\pi :Z \ra X$ such that
$\pi$ is cell-like, and $\dim Z \leq n$.
\end{tm}

We intend to explain this equivalence in Section 2.\\

However, the proof of Theorem \ref{T} in \cite{To} is interesting because it can be done without using Edwards-Walsh complexes, which were used in the original proof of Theorem~\ref{EdWa}. This requires changing the proof of Theorem~3.9 from \cite{To}, which will be done in Section 3 of this paper. 

The definition and properties of Edwards-Walsh complexes can be found in \cite{Dr1}, \cite{DW} or \cite{KY}. Using Edwards-Walsh complexes, or CW-complexes built similarly to these, was the standard approach in proving resolution theorems, for example in \cite{Wa}, \cite{Dr1} and \cite{Le}. But these complexes can become fairly complicated, which also complicates the algebraic topology machinery appearing in proofs using them. The proof of Theorem~\ref{T}, after the adjustment of proof of Theorem~3.9 from \cite{To}, does not use Edwards-Walsh complexes -- instead, it has a more involved point set topological part.
Therefore the Edwards-Walsh resolution theorem can be proven without using Edwards-Walsh complexes.

\section{The equivalence of the two theorems}

We will use the following theorem by A.~Dranishnikov, which can be found in \cite{Dr1} as Theorem 11.4, or in \cite{Dr2} as Theorem 9:

\begin{tm}\label{Dran1}
For any simple \emph{CW}-complex $M$ and any finite dimensional
compactum $X$, the following are equivalent:
\begin{enumerate}
\item $X\tau M$;
\item $X\tau SP^{\infty}M$;
\item $\dim_{H_i(M)} X \leq i$ for all $i\in \N$;
\item $\dim_{\pi_i(M)} X \leq i$ for all $i\in \N$.
\end{enumerate}
\end{tm}

A space $M$ is called \emph{simple} if the action of the fundamental group $\pi_1(M)$ on all homotopy groups is trivial. In particular, this implies that $\pi_1(M)$ is abelian.
Also, $SP^\infty M$ is the infinite symmetric product of $M$, and for a CW-complex $M$, $SP^\infty M$ is homotopy equivalent to the weak cartesian product of Eilenberg-MacLane complexes $K(H_i(M),i)$, for all $i\in \N$.

In fact, Theorem 6 from \cite{Dr2} states that if $X$ is a compact metrizable space, and $M$ is any CW-complex, then $X \tau M$ implies $X\tau SP^\infty M$. Moreover, since $SP^\infty M$ is homotopy equivalent to the weak  product of Eilenberg-MacLane complexes $K(H_i(M ),i)$, then $X\tau SP^\infty M$ implies $X \tau K(H_i(M ),i)$, for all $i\in \N$. This means that the implications 
 (1) $\Rightarrow$ (2) $\Rightarrow$
(3) from Theorem~\ref{Dran1} are true for any compact metrizable space $X$, and not just for finite dimensional ones, as well as for any CW-complex $M$. So we can restate a part of the statement of  Theorem~\ref{Dran1} in the form we will need:

\begin{tm}\label{Dran2}
For any \emph{CW}-complex $M$ and any 
compact metrizable space $X$, we have
$X\tau M$ $\Rightarrow$ $X\tau SP^{\infty}M$ $\Rightarrow$
$\dim_{H_i(M)} X \leq i$ for all $i\in \N$.
\end{tm}

Now $X$ from Theorem~\ref{T} has property $X\tau K$, where $K$ is a connected CW-complex with $\pi_n(K)\cong G$, $\pi_k(K)\cong 0$ for $0\leq k< n$, and $n\in \N$.
 By Hurewicz Theorem, if $n=1$, since $G$ is abelian we get $H_1(K)\cong\pi_1(K)$, and if $n\geq 2$ then $H_n(K)\cong\pi_n(K)$. Therefore,  by Theorem~\ref{Dran2}, $X \tau K$ implies $\dim_{H_n(K)} X \leq n$, i.e.,  $\dim_G X\leq n$.

By Bockstein Theorem and basic properties of Bockstein basis, as explained in Lemma~2.4 from \cite{To}, $P_G=\P$ implies that $\dim_G X=\dim_\Z X$. Now use the Edwards-Walsh resolution theorem to produce a compact metrizable space $Z$ with $\dim Z\leq n$, and a cell-like map $\pi: Z\ra X$. Since $\dim_A Z\leq \dim Z$ for any abelian group $A$, using $A=H_n(K)=G$ as well as other properties of $K$, and the fact that $Z$ is finite dimensional, Lemma~3.10 from \cite{To} shows $Z\tau K$.

\section{How to avoid using Edwards-Walsh complexes}

In the proof of Theorem~\ref{T} in \cite{To}, the following theorem is used -- it appears in \cite{To} as Theorem~3.9. This theorem is a known result, presented in a particular form that was adjusted to fit the needs of the proof of Theorem~\ref{T}. This is why its proof was presented in \cite{To}.

\begin{tm}[A variant of Edwards' Theorem] \label{Ed}
Let $n\in \N$ and let $Y$ be a
compact metrizable space such that $Y=\lim \ (\vert L_i\vert,
f_i^{i+1})$, where $\vert L_i\vert$ are compact polyhedra with
$\dim L_i \leq n+1$, and $f_i^{i+1}$ are surjections. Then
$\dim_\Z Y \leq n$ implies that there exists an $s\in \N$, $s >1$,
and there exists a map $g_1^s : \vert L_s\vert \to \vert
L_1^{(n)}\vert$ which is an $L_1$-modification of $f_1^s$.
\end{tm}

The proof of this theorem in \cite{To} had two parts, the first part for $n\geq 2$ and the second for $n=1$. In the first part of the proof,  Edwards-Walsh complexes were used.
The proof is still correct, but it turns out that there was no need to use Edwards-Walsh complexes. In fact, the entire proof can be simplified, and done for any $n\in\N$ as it was done for the case when $n=1$.
Theorem~\ref{Ed} was the only place in \cite{To} where Edwards-Walsh complexes were used, so the main result of \cite{To} can be proven without ever using them. Consequently, the Edwards-Walsh resolution theorem can be proven without using Edwards-Walsh complexes.

The goal of this section is to give a simplified proof for Theorem~\ref{Ed}.
Here is a reminder of some facts from the original paper that are used in the new proof.

\vspace{2mm}

First of all, recall that a map $g:X \ra \vert K\vert$ between a space $X$ and a simplicial complex $K$ is called a $K$-\emph{modification} of $f$ if whenever $x\in X$ and $f(x)\in \sa$, for some $\sa \in K$, then $g(x)\in \sa$. This is equivalent to the following: whenever $x\in X$ and $f(x)\in \overset{\circ}\sa$, for some $\sa \in K$, then $g(x)\in \sa$.

\vspace{2mm}

In the course of the simplified proof of Theorem~\ref{Ed}, we will need the notion of \textit{resolution in the sense of inverse sequences}. This usage of the word resolution is completely different from the notion from the title of this paper. The definition can be found in \cite{MS} for the more general case of inverse systems. We will give the definition for inverse sequences only.

Let $X$ be a topological space. A \emph{resolution} of $X$ \emph{in the sense of inverse sequences} consists of an inverse sequence of topological spaces $\mathbf{X}= (X_i,p_i^{i+1})$ and a family of maps $(p_i:X \ra X_i)$ with the following two properties:
\begin{enumerate}
\item[(R1)] Let $P$ be an ANR, $\scv$ an open cover of $P$ and $h:X\ra P$ a map. Then there is an index $s\in \N$ and a map $f:X_s \ra P$ such that the maps $f\circ p_s$ and $h$ are $\scv$-close.
\item[(R2)] Let $P$ be an ANR and $\scv$ an open cover of $P$. There exists an open cover $\scv '$ of $P$ with the following property: if $s\in\N$ and $f,f':X_s \ra P$ are maps such that the maps $f\circ p_s$ and $f'\circ p_s$ are $\scv '$-close, then there exists an $s'\geq s$ such that the maps $f\circ p_s^{s'}$ and $f'\circ p_s^{s'}$ are $\scv$-close.
\end{enumerate}

By Theorem I.6.1.1 from \cite{MS}, if all $X_i$ in $\mathbf{X}$ are compact Hausdorff spaces, then $\mathbf{X}= (X_i,p_i^{i+1})$ with its usual projection maps $(p_i:\lim \mathbf{X} \ra X_i)$ is a resolution of $\lim \mathbf{X}$ in the sense of inverse sequences.
Moreover, since every compact metrizable space $X$ is the inverse limit of an inverse sequence of compact polyhedra $\mathbf{X}= (P_i,p_i^{i+1})$ (see Corollary I.5.2.4 of \cite{MS}), this inverse sequence $\mathbf{X}$ will have the property (R1) mentioned above, and we will refer to this property as the \emph{resolution property} (R1) \emph{in the sense of inverse sequences}.

\vspace{2mm}

We will also use stability theory, about which more details can be found in \S VI.1 of \cite{HW}. Namely, we will use the consequences of Theorem VI.1. from \cite{HW}: if $X$ is a separable metrizable space with $\dim X \leq n$,  then for any map $f: X\ra I^{n+1}$, all values of $f$ are unstable. A point $y\in f(X)$ is called an  \emph{unstable value} of $f$ if for every $\da > 0$ there exists a map $g:X\ra I^{n+1}$ such that:
\begin{enumerate}
\item $d (f(x),g(x))<\da$ \ for every $x\in X$, and
\item$g(X)\subset I^{n+1}\setminus\{y\}$.
\end{enumerate}
Moreover, this map $g$ can be chosen so that $g=f$ on the complement of $f^{-1}(U)$, where $U$ is an arbitrary open neighborhood of $y$, and so that $g$ is homotopic to $f$ (see Corollary I.3.2.1 of \cite{MS}).

\vspace{2mm}

Here is a technical result from \cite{To}, which is stated there as Lemma~3.7 and used in the proof of Theorem~\ref{Ed}.

\begin{lemma}\label{sh}
For any finite simplicial complex $C$, there is a map $r:\vert C\vert \ra \vert C\vert$ and an open cover
$\mathcal{V}=\{V_\sa : \sa \in C\}$ of $\vert C\vert$ such that for all $\sa$, $\tau \in C$:
\begin{enumerate}
\item[(i)] $\overset{\circ}\sa \subset V_\sa$,
\item[(ii)] if $\sa\neq \tau$ and $\dim \sa = \dim \tau$, $V_\sa$ and $V_\tau$ are disjoint,
\item[(iii)] if $y\in \overset{\circ}\tau$,
 $\dim \sa \geq \dim \tau$ and $\sa \neq \tau$, then  $y \notin V_\sa$,
\item[(iv)]if $y \in \overset{\circ}\tau \cap V_\sa$, where  $\dim \sa <\dim \tau$, then $\sa$ is a face of $\tau$, and
\item[(v)] $r(V_\sa)\subset \sa$.
\end{enumerate}
\end{lemma}


\noindent\textit{Simplified proof of Theorem~\ref{Ed}}:
Since $Y=\lim (\vert L_i\vert, f_i^{i+1})$, where $\vert L_i\vert$ are compact polyhedra with $\dim L_i\leq n+1$, we get that $\dim Y \leq n+1$. According to Aleksandrov's theorem (\cite{Al}), $\dim Y$ being finite means $\dim_{\Z} Y = \dim Y$.
Therefore, assuming $\dim_{\Z} Y \leq n$ really means that $\dim Y \leq n$, too.

Thus we can prove the theorem without using Edwards-Walsh complexes, but instead using the resolution property (R1) in the sense of inverse sequences.

We can construct a map  $g_1:Y \ra \vert L_1^{(n)}\vert$ that equals $f_1$ on $f_1^{-1}(|L_1^{(n)}|)$.
This can be done as follows.
Let $\sa$ be an $(n+1)$-simplex of $L_1$ and $w\in
\overset{\circ}\sa$. Since $\dim \sa =n+1$ and $\dim Y \leq n$, the point $w$ is an
unstable value for $f_1$ ($f_1$ is surjective, since all our bonding maps $f_i^{i+1}$ are surjective). Therefore we can find a map $g_{1,\sa}:Y\ra \vert
L_1\vert$ which agrees with $f_1$ on
$Y\setminus (f_1^{-1}(\overset{\circ}\sa))$, and $w \notin g_{1,\sa}(Y)$.
Then choose a map $r_\sa :\vert L_1\vert \ra \vert L_1\vert$ such that
$r_\sigma$ is the identity on $|L_1|\setminus \overset{\circ}\sa$ and $r_\sigma(g_{1,\sigma}(Y))\cap \overset{\circ}\sa=\emptyset$. Finally, replace $f_1$ by $r_\sigma\circ g_{1,\sigma}: Y \ra \vert L_1\vert \setminus \overset{\circ}\sa$.

Continue the process with one $(n+1)$-simplex at a time. Since $L_1$
is finite, in finitely many steps we will reach the needed map
$g_1:Y \ra \vert L_1^{(n)}\vert$. Note that from the construction of
$g_1$, we get
\begin{enumerate}
\item[(I)] $g_1\vert_{f_1^{-1}(\vert
L_1^{(n)}\vert)}=f_1\vert_{f_1^{-1}(\vert L_1^{(n)}\vert)}$, and
for every $(n+1)$-simplex $\sa$ of $L_1$, $\ g_1(f_1^{-1}(\sa))\subset
\partial \sa$.
\end{enumerate}
\begin{displaymath}
\xymatrix{
 \vert L^{(n)}_1 \vert \ar@{_{(}->}[d] & & & &\\
\vert L_1\vert &  & \vert L_s\vert \ar[ll]^{f_1^s} \ar@{-->}[ull]_{\hspace{-3mm}\widehat{g}_1^s} \ar@/^/@{.>}[ull]|{g_1^s} & ...\ar[l]  & Y \ar@/^/[ll]^{\qquad f_s}  \ar@/_/[ullll]_{g_1} \ar@/^2pc/[llll]^{\qquad f_1}\\
&&&&&
}
\end{displaymath}

Let us choose an open cover $\mathcal{V}$ of $\vert L_1^{(n)}\vert$ by applying Lemma~\ref{sh} to $C=L_1^{(n)}$.
Now we can use  resolution property (R1) in the sense of inverse sequences: there is an index $s>1$
and a map $\widehat{g}_1^s:\vert L_s \vert \ra \vert
L_1^{(n)}\vert$ such that $\widehat{g}_1^s\circ f_s$ and $g_1$ are
$\mathcal{V}$-close. Define $g_1^s:=r\circ \widehat{g}_1^s
:\vert L_s \vert \ra \vert L_1^{(n)}\vert$, where $r: \vert L_1^{(n)}\vert \ra \vert L_1^{(n)}\vert$ is the map from Lemma~\ref{sh}.

Notice that for any $y\in Y$, if $g_1(y)\in \overset{\circ}\tau$
for some $\tau \in L_1^{(n)}$, then $g_1(y)\in
V_\tau$, and possibly also $g_1(y) \in V_{\ga_j}$, where $\ga_j$ are some faces of
 $\tau$ (there can only be finitely many). Then either $\widehat{g}_1^s \circ f_s(y) \in V_\tau$,
or $\widehat{g}_1^s \circ f_s(y) \in V_{\ga_j}$, for some $\ga_j$. In any case, $r\circ
\widehat{g}_1^s \circ f_s(y) \in \tau$. Hence,
\begin{enumerate}
\item[(II)] for any $y\in Y$, $g_1(y)\in \overset{\circ}\tau$ for
some $\tau \in L_1^{(n)}$ implies that $g_1^s(f_s(y)) \in \tau$.
\end{enumerate}

Finally, for any $z \in \vert L_s\vert$, $f_s$ is surjective
implies that there is a $y \in Y$ such that $f_s(y)=z$. Then
$f_1^s(z)=f_1^s(f_s(y))=f_1(y)$. Now $f_1^s(z)$ is either in
$\overset{\circ}\sa$ for some $(n+1)$-simplex $\sa$ in $L_1$, or in
$\overset{\circ}\tau$ for some $\tau \in L_1^{(n)}$.

If $f_1^s(z)\in \overset{\circ}\sa$, that is $f_1(y) \in
\overset{\circ}\sa$ for some $(n+1)$-simplex $\sa$, by (I) we get
that $g_1(y)\in
\partial \sa$. Then by (II), $g_1^s(f_s(y)) \in \partial \sa$,
i.e., $g_1^s(z) \in \sa$.

If $f_1^s(z)=f_1(y)\in \overset{\circ}\tau$ for some $\tau \in
L_1^{(n)}$, then (I) implies that $g_1(y)=f_1(y) \in
\overset{\circ}\tau$, so by (II), $g_1^s(f_s(y)) \in \tau$, i.e.,
$g_1^s(z)\in \tau$.

Therefore, $g_1^s$ is indeed an $L_1$-modification of
$f_1^s$.\hfill $\square$

\section{A note about the original proof of the Edwards-Walsh resolution theorem}

In the original proof of Theorem~\ref{EdWa} in \cite{Wa}, the following theorem is used.
It is listed there as Theorem~4.2.

\begin{tm}[R.\ Edwards] \label{EdOriginal}
Let $n\in \N$ and let $X$ be a
compact metrizable space such that $X=\lim \ (P_i,
f_i^{i+1})$, where $P_i$ are compact polyhedra.
The space $X$ has cohomological dimension $\dim_\Z X \leq n$ if and only if for each integer $k$ and each $\e >0$ there is an integer $j>k$, and a triangulation $L_k$ of $P_k$ such that for any triangulation $L_j$ of $P_j$ there is a map $g_k^j: \vert L_j^{(n+1)}\vert \ra \vert L_k^{(n)}\vert$ which is $\e$-close to the restriction of $f_k^j$.
\end{tm}

There were no additional assumptions made about dimension of polyhedra $P_i$, so in the proof of this theorem  in \cite{Wa}, the usage of Edwards-Walsh complexes is indispensable. Therefore, the usage of Edwards-Walsh complexes was necessary in the original proof of Theorem~\ref{EdWa} in \cite{Wa}.

Theorem~\ref{Ed} was modeled on Theorem~\ref{EdOriginal}, but with the additional assumption about dimension of polyhedra $\dim \vert L_i \vert\leq n+1$. This assumption, together with $\dim_\Z Y \leq n$ implies that $\dim Y \leq n$. Therefore the usage of Edwards-Walsh complexes in its proof can be avoided altogether. In fact, Theorem~\ref{Ed} becomes analogous to  Theorem~4.1 from \cite{Wa} -- a weaker version of Edwards' Theorem:

\begin{tm}
Let $n\in \N$ and let $X$ be a
compact metrizable space such that $X=\lim \ (P_i,
f_i^{i+1})$, where $P_i$ are compact polyhedra. The space $X$ has $\dim X \leq n$ if and only if for each integer $k$ and and each $\e >0$ there is an integer $j>k$, a triangulation $L_k$ of $P_k$, and a map $g_k^j: P_j \ra \vert L_k^{(n)}\vert$ which is $\e$-close to  $f_k^j$.
\end{tm}

\noindent\textbf{Acknowledgements}. The author cordially thanks the anonymous referee for their
valuable comments which lead to significant improvement of this paper.

\end{document}